\documentclass[smallextended,numbook,runningheads]{svjour3}
\usepackage{amsmath}
\smartqed  % flush right qed marks, e.g. at end of proof

\usepackage{amsfonts,amssymb}
\usepackage{epsfig}
\usepackage{booktabs}
\usepackage{xcolor}

\journalname{}
\date{ \phantom{b} \vspace{45mm}\phantom{e}}
\def\makeheadbox{\hspace{-4mm}Version of 24 May 2022}%, file: adiabatic.tex}

\def\real{{\mathbb R}}

\def\d{{\mathrm d}}
\def\e{{\mathrm e}}

\def\iu{\mathrm{i}}
\def\eps{\varepsilon}

\def\forall{\qquad\hbox{ for all }\ }

\newdimen\GGGlength
\newdimen\GGGheight
\newbox\GGGbox
\def\GGGput[#1,#2](#3,#4)#5{%
  \setbox\GGGbox\vbox{\hbox{#5}\kern0pt}%
  \GGGlength\wd\GGGbox%
  \divide\GGGlength by100 \multiply\GGGlength by#1%
  \GGGheight\ht\GGGbox%
  \divide\GGGheight by100 \multiply\GGGheight by#2%
  \put(#3,#4){\kern-\GGGlength\raise-\GGGheight\box\GGGbox}}

\begin{document}

\title{On a large-stepsize %modified Boris 
integrator for charged-particle dynamics}

\titlerunning{Large-stepsize integrator for charged-particle dynamics}

\author{Christian Lubich$^1$, Yanyan~Shi$^{1}$}
\authorrunning{Ch.\ Lubich, Y.\ Shi}

\institute{
%$^1$~Dept.\ de Math{\'e}matiques, Univ.\ de Gen{\`e}ve,
%CH-1211 Gen{\`e}ve 24, Switzerland.\\
%\phantom{$^1$~}\email{Ernst.Hairer@unige.ch}\\
$^1$~Mathematisches Institut, Univ.\ T\"ubingen, D-72076 T\"ubingen, Germany.\\
\phantom{$^1$~}\email{\{Lubich, Shi\}@na.uni-tuebingen.de}\\
%$^3$~LSEC, Academy of Mathematics and Systems Science, Chinese Academy of Sciences,\\ 
%\phantom{$^3$~}Beijing 100190, China; University of Chinese Academy of Sciences, Beijing 100049, China.\\
%\phantom{$^3$~}\email{shiyanyan1995@lsec.cc.ac.cn}
}

\date{ }

\maketitle

 \begin{abstract} Xiao and Qin [Computer Physics Comm., 265:107981, 2021] recently proposed a remarkably simple modification of the Boris algorithm to compute the guiding centre of the highly oscillatory motion of a charged particle with step sizes that are much larger than the period of gyrorotations. 
% to use the Boris integrator with an appropriately modified electric field and modified initial velocity to compute the guiding centre of the motion of a charged particle, using step sizes that are much larger than the period of gyrorotations. 
 They gave strong numerical evidence %for this approach 
 but no error analysis. This paper provides an analysis of the large-stepsize modified Boris method in a setting that has a strong non-uniform magnetic field and moderately bounded velocities, considered over a fixed finite time interval. The error analysis is based on comparing the modulated Fourier expansions of the exact and numerical solutions, for which the differential equations of the dominant terms are derived explicitly.
 Numerical experiments illustrate and complement the theoretical results.
\bigskip

\end{abstract}

\section{Introduction}
Integrating the equations of motion of charged particles is a fundamental computational task in particle methods of plasma physics, e.g.~\cite{birdsall05ppv}. The standard numerical integrator for these computations is the Boris algorithm \cite{boris70rps}, which has the charm of simplicity and remarkable conservation properties \cite{qin13wib,hairer18ebo}. In the practically important situation of a strong magnetic field and moderate velocities, which will be considered in this paper, the particle trajectories show fast gyrorotations of small radius around a guiding centre \cite{northrop63tam}. 

To approximate the guiding centre motion, one approach --- not considered here --- is to integrate numerically the known but structurally complicated differential equations for the approximate guiding centre by a suitable numerical method \cite{ellison18dvi}. 

In a different and arguably more efficient approach, a Boris-type integrator with appropriate modifications is applied to the original equations of motion of the charged particle with {\it large step sizes} that do not resolve the high-frequency oscillations. As the standard Boris algorithm used with large step sizes is known to produce numerical solutions with unphysically large gyroradius \cite{parker91nei,ricketson20aec}, modifications to it are necessary. In the case of a near-uniform strong magnetic field, it suffices to filter out the normal component of the initial velocity \cite{hairer22lsi}, but the mere modification of initial values is not sufficient in a strongly non-uniform magnetic field. Xiao \& Qin \cite{xiao21smc} recently proposed to additionally modify the electric field in a non-obvious way when using the Boris algorithm with large step sizes and showed striking numerical results, but no error analysis was given. It was then found that a very similar numerical method, with the same extra force term, was already proposed by Vu \& Brackbill \cite{vu95ans} (Method III) in 1995, motivated by Parker \& Birdsall~\cite{parker91nei}  on the large-stepsize behaviour of the Boris method; see in particular formula (10) in \cite{parker91nei}, based on the equations of guiding-centre motion as given by Northrop~\cite{northrop63tam}, which are at the origin of the extra force term added to the schemes in \cite{vu95ans} and \cite{xiao21smc}. A very interesting approach to understanding such methods in terms of slow manifolds was recently given by Burby \& Klotz \cite{burby20smr} and Burby \& Hirvijoki \cite{burby21nso}, but these papers deal with the exact flow on and near the slow manifold and do not clarify the behaviour of the numerical method for large step sizes.

{\it The objective of the present paper is to give a rigorous analysis of the modified Boris algorithm of Xiao \& Qin \cite{xiao21smc} for approximating the guiding centre motion of a charged particle in a strong non-uniform magnetic field taking large time steps that cover many periods of gyrorotation.}

In Section~\ref{sec:main} we formulate the general setting of charged-particle motion in a strongly non-uniform strong magnetic field, describe the Boris algorithm and its modification, and state our main result, Theorem~\ref{thm:main}, which yields second-order error bounds for the position and velocity of the guiding centre when approximated with the modified Boris method with large step sizes whose square exceeds the inverse of the strength of the magnetic field. In Section~\ref{sec:num} we present results of numerical experiments that illustrate and complement the theory.
In Section~\ref{sec:mfe} we give modulated Fourier expansions of both the exact and the numerical solution. Their comparison yields the proof of Theorem~\ref{thm:main}.

\section{Large-stepsize modified Boris method and its error bound}
\label{sec:main}

\subsection{Setting}
\label{subsec:setting}

The motion of a charged particle (of unit mass and charge) in a magnetic and electric field is governed by the differential equation
\begin{equation}\label{ode}
\ddot x =  \dot x \times B(x) + E(x),
\end{equation}
where $x(t)\in\mathbb R^3$ is the position at time $t$, $v(t)=\dot x(t)$ is the velocity, $B$ is the magnetic field and $E$ is the electric field. Here, $B(x) = \nabla \times A(x)$ with a vector potential $A(x)\in\real^3$ and 
$E(x) = - \nabla \phi(x)$ with a scalar potential $\phi(x)\in\real$, which we assume to be bounded from below. 
We are interested in the case of a strong magnetic field
\begin{equation}\label{B-eps}
B(x) = B_\eps(x)= \frac 1\eps \, B_1(x),  \quad\ 0<\eps\ll 1,
\end{equation}
where $B_1$ is smooth and independent of the small parameter $\eps$, with $|B_1(x)|\ge 1$ for all $x$. The motion \eqref{ode} and its approximation are to be studied over time intervals $t\in [0,T]$ with fixed $T$ independent of $\eps$, for initial values $(x(0),\dot x(0))$ that are bounded independently of $\eps$: for some constants $M_0,M_1$,
\begin{equation} \label{xv-bounds}
|x(0)| \le M_0, \quad\ |\dot x(0)| \le M_1.
\end{equation}
Under these conditions, it is known that the magnetic moment
$$
\mu(x,v) = \frac{1}{2}\frac{|v\times B(x)|^2}{|B(x)|^3},
$$
which is of size $O(\eps)$ under our assumptions, is an adiabatic invariant \cite{kruskal58tgo,northrop63tam}: $\mu(x(t),\dot x(t))$ is conserved up to $O(\eps^2)$ over very long times $t\le \eps^{-N}$ with arbitrary $N>1$ \cite{benettin94aia,hairer20lta}. Here we consider \eqref{ode} only over fixed times $T$ that are independent of $\eps$.

\subsection{Modified Boris method of Xiao \& Qin \cite{xiao21smc} }
The integrator for charge-particle dynamics \eqref{ode} proposed in \cite{xiao21smc} is a remarkably simple but nontrivial modification of the Boris algorithm, with the objective to approximate the guiding centre of the particle motion with large step-sizes $h\gg \eps$ without resolving the gyrorotations. In its two-step formulation the method computes the new position $x^{n+1}$ as an approximation at time $t_{n+1}=(n+1)h$ via
\begin{equation}\label{mboris}
\frac{x^{n+1}-2x^n+x^{n-1}}{h^2}=v^n \times B(x^n) + E(x^n)- \mu^0\, \nabla |B|(x^n)
\end{equation}
with the initial magnetic moment $\mu^0=\mu(x(0),\dot x(0))$ and the symmetric finite difference velocity approximation
\begin{equation} \label{vn}
v^n = \frac{x^{n+1}-x^{n-1}}{2h}.
\end{equation}
This differs from the original Boris method only in the addition of the extra term $- \mu^0\, \nabla |B|(x^n)$, which also appears in the differential equations for the guiding centre; see  \cite{northrop63tam} and, e.g.,  Theorem~\ref{thm:mfe} below.
The modifed Boris method starts from modified initial values
\begin{equation}\label{mod-init}
x^0 = x(0), \quad v^0 = P_\parallel(x^0) \,\dot x(0),
\end{equation}
where $P_\parallel(x^0)$ is the orthogonal projection onto the span of $B(x^0)$. With $P_\perp(x^0)=I-P_\parallel(x^0)$, we note that
$$
v^0_\perp = P_\perp(x^0) v^0 = 0,
$$
i.e., the perpendicular component of the initial velocity has been filtered out.

The scheme \eqref{mboris} is identical to the standard Boris integrator for the modified force field 
$$
E_\mathrm{mod}(x) = E(x)- \mu^0 \,\nabla |B|(x) =-\nabla (\phi + \mu^0 |B|)(x).
$$
The actual implementation uses the common one-step formulation of the Boris algorithm \cite{boris70rps}.

\subsection{Large-stepsize error bound}

For the following theorem, which is the main result of this paper,  we need a nondegeneracy condition:
\begin{align}\label{ass-nondeg}
&\text{For $(x,v)$ along the numerical trajectory, the linear maps} \nonumber
\\[1mm]
& \text{$L_{x,v}:P_\perp(x)\mathbb{R}^3 \to P_\perp(x)\mathbb{R}^3, \quad z \mapsto z + \tfrac14 h^2\, P_\perp(x)\bigl(v \times B'(x)z\bigr)$}
\\[1mm]
&\text{have an inverse that is bounded independently of $(x,v)$ and of 
%$h$ and $\eps$ with $h^2/\eps\le C_*$.
} \nonumber
\\[-1mm]
&\text{$h$ and $\eps$ with $h^2/\eps\le C_*$.} \nonumber
%&\text{are invertible with a bound of the inverse that is independent of $h$ and $\eps$.} \nonumber
%\\[-1mm]
%&\text{independent of $h$ and $\eps$.}  \nonumber
\end{align}
This determines an upper bound $C_*$ on the ratio $h^2/\eps$.
We have the following large-stepsize error bound for the modified Boris method.

\begin{theorem}\label{thm:main}
Consider applying the modified Boris method to \eqref{ode}--\eqref{xv-bounds} with modified initial values \eqref{mod-init} over a time interval $0\leq t\leq T$ (with $T$ independent of $\varepsilon$) using a step size $h$ with $h^2\sim \eps$, i.e.,
\[
%h^2\geq \varepsilon.
c_* \eps \le h^2 \le C_* \eps
\]
for some positive constants $c_*$ and $C_*$. Under the nondegeneracy condition \eqref{ass-nondeg},
%Suppose that the numerical solution $x^n$ stays in a compact set $K$ and the starting velocity $v^0$ is bounded independent of $\varepsilon$ and $h$ and its component orthogonal to $B(x(0))$ is zero (or $O(\eps)$).
%Then, 
the errors in position $x$ and parallel velocity $v_\parallel=P_\parallel(x) v$ (where $P_\parallel(x)$ denotes the orthogonal projection onto the span of $B(x)$) at time $t_n=nh\leq T$ are bounded by
\[
|x^n-x(t_n)|\leq Ch^2, \quad |v_\parallel^n-v_\parallel(t_n)|\leq Ch^2 ,\quad\ |v^n_\perp| \le Ch^2,
\]
where $C$ is independent of $\varepsilon$, $h$ and $n$ with $nh\leq T$ (but depends on T, on bounds of derivatives of $B_1$ and $E$, and on $c_*$ and $C_*$).
\end{theorem}

Since $x(t)$ and $v_\parallel(t)$ are $O(\eps)$ close to the guiding centre at time $t$  and its velocity, respectively, and since $\eps \sim h^2$ by assumption, Theorem~\ref{thm:main} yields that the modified Boris method approximates the guiding centre motion with $O(h^2)$ accuracy for step sizes $h$ that are much larger than the gyroperiod $2\pi/|B(x)| \sim \eps$.

The proof of this theorem will be given in Section~\ref{sec:mfe}.

\begin{remark}
It is of interest to understand how the error bound changes when $\eps \ll h^2$ or $h \gg \eps \gg h^2$. 
%The same proof with additional attention to some error terms then yields an error bound of  $O(\rho^2 h^2)$. 
An $O(h^2)$ error bound still holds true for 
the less restrictive stepsize condition
$$
c_* \eps \le h^2 \le C_* \eps^\alpha \qquad\text{for } 0\le \alpha <1
$$
for less strongly varying magnetic fields
$$
B(x)=\frac1\eps\, B_1(\eps^{1-\alpha} x).
%B'(y) = O(\eps^{-\alpha}).   
$$
This can still be obtained with essentially the same proof, but we will not work out the lengthy yet conceptually straightforward details.

In the situation of $h^2 \le \eps \ll h$ given by
$$
c_* \eps^\beta \le h^2 \le C_* \eps \qquad\text{for } 1<\beta<2,
$$
an $O(\eps)$ error bound can be shown without an extra assumption on derivatives of $B$, again with essentially the same proof.
\end{remark}

\section{Numerical experiments}
\label{sec:num}

We show numerical results of the modified Boris method for two examples.

\subsection{Tokamak example from~\cite{xiao21smc}}
We consider the motion of a charged particle in a tokamak geometry without electric field~\cite{xiao21smc}. 
%The magnetic field in the toroidal coordinates $(r,\theta,\xi)$ of $x=(x_1,x_2,x_3)$ given by
%\[
%r=\sqrt{(R-R_0)^2+x_3^2}, \quad \theta=\arctan\left(\frac{x_3}{R-R_0}\right), \quad \xi=\arctan\left(\frac{x_2}{x_1}\right)
%\] 
%with $R=\sqrt{x_1^2+x_2^2}$ is expressed as
%\[
%B(x)=\frac{B_0r}{qR}\e_{\theta}+\frac{B_0R_0}{R}\e_{\xi}.
%\]
%We choose $B_0=1$, $R_0=1$, and $q=2$. Equivalently in 
In Cartesian coordinates, the magnetic field is given as
\[
B(x)=\Bigl(-\frac{2x_2+x_1 x_3}{2R^2},\ \frac{2x_1-x_2x_3}{2R^2}, \ \frac{R-1}{2R} \Bigr)^\top
\quad\text{ with } \ R=\sqrt{x_1^2+x_2^2}\,.
\]
Starting with the initial position $x(0)=(1.05,0,0)^\top$ and the initial velocity $\dot{x}(0)=(2.1\times10^{-3}, 4.3\times10^{-4},0)^\top$, the orbit projected onto the $(R,x_3)$ plane is a banana orbit.  The final time considered is $T=3.75\times10^4$. We note that upon rescaling time as $t \to \eps t$ with $\eps=10^{-3}$, the problem has the scaling of Section~\ref{subsec:setting}.

\begin{figure}[h]
\centerline{\includegraphics[scale=0.5]{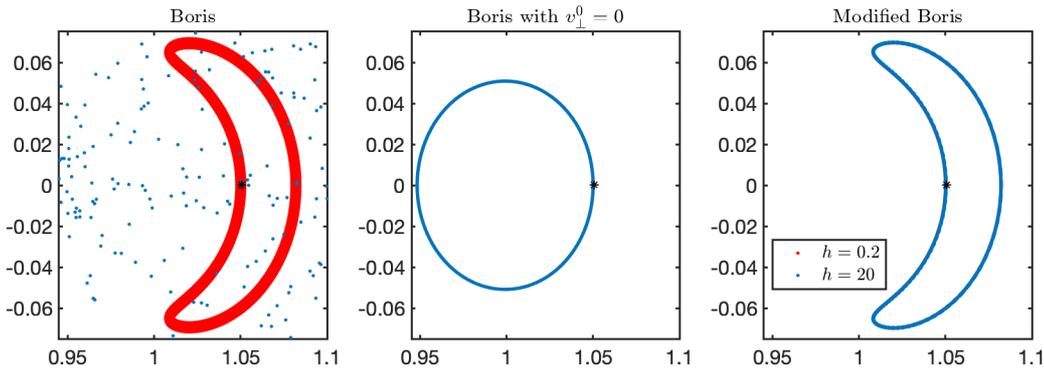}}
\caption{Banana orbits  on the $R-x_3$ plane computed by different methods with $T=3.75\times10^4$, $h=0.2$ (red) and $h=20$ (blue).}\label{fig:orbit}
\end{figure}

Figure~\ref{fig:orbit} shows the trajectories computed by the standard Boris, standard Boris with projected initial velocity and the modified Boris algorithm.  Two step sizes  $h=0.2$ and $h=20$ are chosen. It is observed that when $h=0.2$, the standard Boris shows the correct result while the gyroradius gets larger with larger step size. For $h=20$ the numerical result is completely wrong. If we use standard Boris with $v^0_\perp=0$, the gyroradius is always small, but the trajectory is not correct. After adding the $\mu\nabla|B|$  term,  the method shows  correct results even for the large step-size $h=20$.

\subsection{Order of accuracy}
We test the numerical accuracy of modified Boris algorithm with large time step by applying the scheme to the example in~\cite{hairer20lta}. We have the electric field
\[
E(x)=-\nabla \phi(x) \quad\text{with the potential }\
\phi(x)=x_1^3-x_2^3+\frac{1}{5}x_1^4+x_2^4+x_3^4,
\]
the magnetic field
\[
B(x)=\frac{1}{2\varepsilon}
\begin{pmatrix}
x_2-x_3\\ x_1+x_3\\ x_2-x_1
\end{pmatrix},
\]
and we take initial values 
\[
x(0)=(0.0, 1.0, 0.1)^\top,\quad \dot{x}(0)=(0.09,0.55,0.3)^\top.
\]

\begin{figure}[h]
\centerline{\includegraphics[scale=0.5]{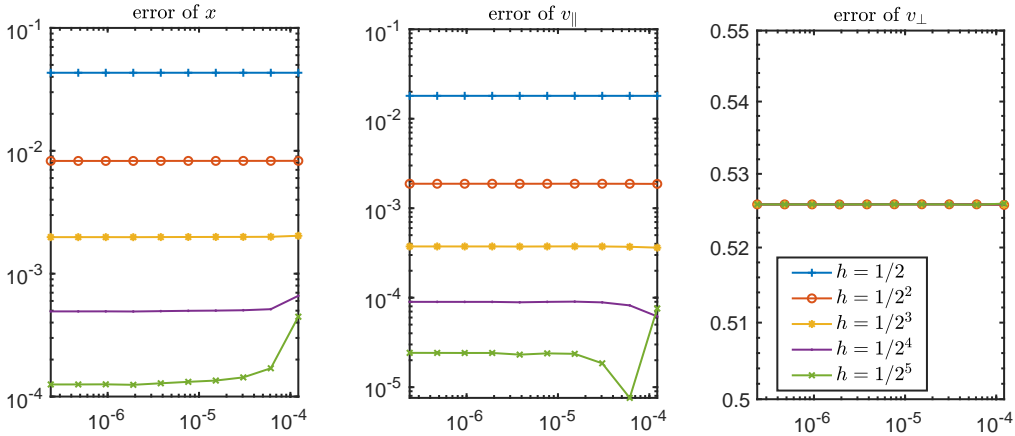}}
\caption{Global error vs. $\varepsilon$ ($\varepsilon=1/2^j,j=13,\cdots 22$) with different $h$ for the modified Boris algorithm.}\label{fig:order}
\end{figure}

Figure~\ref{fig:order} shows the absolute errors in $x$, $v_\parallel$ and $v_\perp$ at time $t=1$ versus $\varepsilon$ for various $h$. It is observed that the errors in $x$ and $v_\parallel$ tend to a constant error level proportional to $h^2$, which is in accordance with our theoretical result in Theorem~\ref{thm:main}.

\section{Modulated Fourier expansions and proof of Theorem~\ref{thm:main}}~\label{sec:mfe}
Theorem~\ref{thm:main} will be proved by comparing the modulated Fourier expansions of the exact and numerical solutions.
Modulated Fourier expansions have previously been used in the analysis of numerical methods for oscillatory differential equations, see \cite{hairer00lec,hairer06gni} and numerous further papers, and lately also for charged-particle dynamics in a strong magnetic field \cite{hairer20lta,hairer22lsi,hairer20afb,wang21eep,wang21eeo}. Incidentally, modulated Fourier expansions (though not under this name) were used for studying the gyration of charged particles by Kruskal \cite{kruskal58tgo} as early as 1958.

The analysis given here builds on that of \cite{hairer20lta} for the exact solution in the situation of a strong non-uniform magnetic field and on that of \cite{hairer22lsi} for the Boris method with large step sizes in the situation of a mildly non-uniform strong magnetic field.

In this section we give the modulated Fourier expansions of the exact solution (Theorem~\ref{thm:mfe}) and of the numerical solution (Theorem~\ref{thm:mfe-boris}), including explicit expressions for the differential equations of the dominant modulation functions. The proof of Theorem~\ref{thm:main} is then readily obtained by a comparison of Theorems~\ref{thm:mfe} and~\ref{thm:mfe-boris}.

\subsection{Modulated Fourier expansion of the exact solution}
We write the solution of \eqref{ode} as
\[
x(t)\approx\sum_{k\in\mathbb Z}z^k(t)\,\e^{{\iu}k\varphi(t)/\eps}, \qquad 0\le t \le T,
\]
with piecewise smooth modulation functions $z^k$ and phase function $\varphi$ for which all time derivatives are bounded independently of $\varepsilon$, except at the discontinuities of $z^k$ and $\dot \varphi$ at integral multiples $t=n\eps$ with jumps of size $O(\eps^N)$, for an arbitrarily chosen integer $N>1$.
The phase satisfies  
$\dot{\varphi}(t)/\eps =|B(z^0(t))|$ 
and $z^0(t)$ is the guiding centre at time $t$ (defined up to $O(\eps^N)$).

Following \cite{hairer20lta}, we diagonalize the linear map $v \mapsto v\times B(x)$, which has eigenvalues $\lambda_1=\iu|B(x)|$, $\lambda_0=0$ and $\lambda_{-1}=-\iu|B(x)|$. The corresponding normalized eigenvectors are denoted by $\nu_1(x)$, $\nu_0(x)$, and $\nu_{-1}(x)$. We let $P_j(x)=\nu_j(x)\nu_j(x)^*$ be the orthogonal projections onto the eigenspaces, where we note that $P_\parallel(x)=P_0(x)$ and $P_\perp(x)=I-P_\parallel(x)= P_1(x)+P_{-1}(x)$. We write the coefficient functions in the basis $\nu_j(z^0(t))$,
\[
z^k=\sum_{j=-1}^1 z^k_j, \qquad z^k_j(t)=P_j(z^0(t))z^k(t), \qquad k\ne 0,
\]
whereas for $k=0$ we decompose
\[
z^0 = c^0 + \sum_{j=-1}^1  z^0_j, \qquad z^0_j(t) = P_j(z^0(t))(z^0(t) -c^0(t)),
\]
where $c^0(t)$ is a piecewise constant function with a finite number of jumps independent of $\eps$, chosen arbitrarily such that $|x(t) - c^0(t)|$ remains distinctly smaller than the inverse of a bound of the derivative of $P_0$ in a neighbourhood of the solution.

The following result is based on Theorem 4.1 of \cite{hairer20lta}, where the existence of the modulated Fourier expansion of solutions of \eqref{ode}--\eqref{B-eps} was established together with bounds of the modulation functions and of the remainder term. However, the differential equations for the dominant modulation functions $z^0_0$, $z^0_{\pm 1}$ and $z^{\pm 1}_{\pm 1}$ and their initial values were not stated explicitly. As these will be needed in the following and are also of independent interest, they are given here.

\begin{theorem}\label{thm:mfe}
Let $x(t)$ be a solution of \eqref{ode}--\eqref{B-eps}  with an initial velocity bounded independently of $\eps$ $(|\dot x(0)|\le M)$.
%, which stays in a compact set $K$ for $0\leq t\leq T$ with $T_{\eps}=O(\eps)$. 
For an arbitrary truncation index $N\geq 1$ we then have an expansion
\[
x(t)=\sum_{|k|\leq N-1}z^k(t)\,\mathrm e^{\mathrm i k \varphi(t)/\varepsilon}+R_N(t),  \qquad 0\le t \le T,
\]
where the phase function satisfies $\dot{\varphi}(t)=|B_1(z^0(t))|$ (recall that $B(x)=B_1(x)/\eps$ with $B_1$ independent of $\eps$), and we fix $\varphi(0)=0$.
%and {\blue $z^k(t)=\sum_m z^{k,[m]}(t)\mathbbm 1_m(t)$ where $\mathbbm 1_m(t)$ is the indicator function defined as
%\[
%\mathbbm 1_m(t)=
%\left\{
%\begin{aligned}
%&1 \quad t\in[m\eps,(m+1)\eps]\\
%&0 \quad \text{other}
%\end{aligned}
%\right..
%\]}
%The following properties hold:
\begin{itemize}
%\item [(a)] 
%\[
%z^{k,[m]}((m+1)\varepsilon)=z^{k,[m+1]}((m+1)\varepsilon)+O(\varepsilon^{N+2}).
%\]
\item [(a)] 
The coefficient functions $z^k(t)$ are piecewise continuous with jumps of size $O(\eps^{N})$ at integral multiples of $\eps$ and are smooth elsewhere. Together with their derivatives (up to order $N$) they are bounded as
\[
z^k=O(\eps^{|k|}) \quad \forall |k|\leq N-1
\]
and further satisfy $z^k_j  =O(\eps^2)$ for $ |k|=1, j\neq k$. Moreover, $\dot z^0 \times B_1(z^0)=O(\eps)$.
%\[
%\dot{z}^0\times B(z^0)=O(\eps), \qquad P_j(z^0)z^k=O(\eps^2) \quad \text{for} \ |k|=1, j\neq k.
%\]
The functions $z^k$ are unique up to $O(\eps^{N})$.

%On every time interval $m\eps\leq t\leq(m+1)\eps$, the modulation functions $z^{k,[m]}$ together with their derivatives (up to order $N$) are bounded as
%\[
%z^{k,[m]}=O(\varepsilon^{|k|}) \  \text{for all} \ |k|\leq N+1
%\]
%and further satisfy 
%\[
%\dot{z}^{0,[m]}\times B(z^{0,[m]})=O(\eps), \quad P_j(z^{0,[m]})z^{k,[m]}=O(\eps^2) \  \text{for} \ |k|=1, j\neq k.
%\]
%They are unique up to $O(\varepsilon^{N+2})$ and are chosen to satisfy $z^{-k,[m]}_{-j} = \overline{z^{k,[m]}_j}$.

\item [(b)] 
The remainder term and its derivative are bounded by
\[
R_N(t)=O(\varepsilon^{N}),\quad \dot{R}_N(t)=O(\varepsilon^{N-1}) , \qquad 0\le t \le T.
\]
\item [(c)]  
On each time interval $n\eps \le t < (n+1)\eps\le T$ (for integers $n\ge 0$), the functions $z_0^{0}$, $z_{\pm1}^{0}$, $z_1^{1}$, $z_{-1}^{-1}$ satisfy the following differential equations. Here, all functions $B$, $E$, $P_j$ are evaluated at the guiding centre $z^0(t)$, and we write $\dot P_j = (d/dt) P_j(z^0(t))=P_j'(z^0(t)) \dot z^0(t)$ and analogously $\ddot P_j$. Moreover, $\mu^0=\mu(x(0),\dot x(0)))$ is the magnetic moment. Omitting the ubiquitous argument $t$, we have
%\[
%\begin{aligned}
%P_0(z^{0,[m]})\ddot{z}^{0,[m]}=&P_0(z^{0,[m]})E(z^{0,[m]})+2P_0(z^{0,[m]}) \, {\rm Re}\blue{ \underbrace{\Bigl(\frac{{\rm i}\dot{\varphi}}{\varepsilon}\,z^{1,[m]}\times \frac 1 \eps B'(z^{0,[m]})z^{-1,[m]}\Bigr)}_{=-\mu|\nabla B(z^{0,[m]})|/2}}\\
%&+2P_0(z^{0,[m]}) \, {\rm Re}\Bigl(\dot{z}^{1,[m]}\times \frac 1 \eps B'(z^{0,[m]})z^{-1,[m]}\Bigr)+O(\varepsilon^2),\\
%P_{\pm 1}(z^{0,[m]})\dot{z}^{0,[m]}=&\pm\iu\frac{\eps}{\dot{\varphi}}\dot{P}_{\pm1}(z^0)\dot{z}^0\pm{\rm i}\frac{\varepsilon}{\dot{\varphi}} P_{\pm 1}(z^{0,[m]})(E(z^{0,[m]})-I\nabla|B|)+O(\varepsilon^2),\\[1mm]
%P_{\pm 1}(z^{0,[m]})\dot{z}^{\pm 1,[m]}=&-\frac{\ddot{\varphi}}{\dot{\varphi}}z^{\pm 1,[m]}_{\pm 1}\mp{\rm i}\frac{\varepsilon}{\dot{\varphi}} P_{\pm 1}(z^{0,[m]})\Bigl(\dot{z}^{0,[m]}\times \frac 1 \eps B'(z^{0,[m]})z^{\pm1,[m]}\Bigr)+O(\varepsilon^2).
%\end{aligned}
%\]
\begin{align*}
\ddot{z}^{0}_0&=2\dot{P}_0 \dot{z}^0+\ddot{P}_0 (z^0-c^0)+ P_0\left(E -\mu^0\,\nabla |B |\right)+O(\eps),\\
\dot{z}^0_{\pm1}&=\dot{P}_{\pm1} (z^0-c^0)\pm\frac{\iu}{|B |}\,\dot{P}_{\pm1} \dot{z}^0\ \pm
\frac{\iu}{|B |}\,P_{\pm1} \left(E -\mu^0\,\nabla |B |\right)+O(\eps^2),\\
\dot{z}^{\pm1}_{\pm1}&=\dot{P}_{\pm1} z^{\pm1}_{\pm1}-\frac{(d/dt)|B |}{|B |}z^{\pm 1}_{\pm1}\ \mp 
 \frac{\iu}{|B |}\,P_{\pm1}\left(\dot{z}^0\times B'(z^0) z^{\pm1}_{\pm1}\right)+O(\eps^2).
\end{align*}
All other modulation functions $z_j^{k}$ are given by algebraic expressions depending on $z^{0}$, $\dot{z}_0^{0}$, $z_1^{1}$, $z_{-1}^{-1}$.
\item [(d)]  
Initial values for the differential equations of item (c) are given by
\[
\begin{aligned}
z^{0}(0)=&\ x(0)+ \frac{\dot{x}(0)\times B(x(0))}{|B(x(0))|^2}+O(\varepsilon^2),\\
\dot{z}_0^{0}(0)=&\ P_0\dot{x}(0) +\dot{P}_0(z^{0}(0)-c^0(0))\ + \\
&\ \frac{1}{|B|^2}\,P_0\bigl(P_0\dot{x}(0)\times B' P_\perp \dot{x}(0)\bigr)+
O(\varepsilon^2),\\
z_{\pm 1}^{\pm 1}(0)=&\ \frac{\mp{\mathrm i}}{|B|}\, P_{\pm 1}\dot{x}(0)+O(\varepsilon^2),
\end{aligned}
\]
where $B,B'$ %are evaluated at $x(0)$ 
and $P_j$ are evaluated at the initial guiding centre $z^0(0)$ (up to $O(\eps^2)$).
%\[
%\begin{aligned}
%z^{0,[m]}(m\eps)=&x(m\eps)+\varepsilon \frac{\dot{x}(m\eps)\times B(x(m\eps))}{|B(x(m\eps))|^2}+O(\varepsilon^2),\\
%\dot{z}_0^{0,[m]}(m\eps)=&P_0(z^{0,[m]}(m\eps))\dot{x}(m\eps) +\dot{P}_0(z^{0,[m]}(m\eps))z^{0,[m]}(m\eps)+ \\
%&\frac{\eps}{\dot{\varphi}^2}P_0(x(m\eps))\bigl(P_0\dot{x}(m\eps)\times(B'(x(m\eps))(P_1 \dot{x}(m\eps)+P_{-1}\dot{x}(m\eps)))\bigr)+
%O(\varepsilon^2),\\
%z_{\pm 1}^{\pm 1,[m]}(m\eps)=&\mp{\mathrm i}\frac{\varepsilon}{\dot{\varphi}(m\eps)} P_{\pm 1}\dot{x}(m\eps)+O(\varepsilon^2).
%\end{aligned}
%\]
\end{itemize}
The constants symbolized by the $O$-notation are independent of $\varepsilon$ and $t$ with $0\leq t\leq T$, but depend on $N$, on the velocity bound $M$, on bounds of derivatives of $B$ and $E$ in a neighbourhood of the trajectory $\{ x(t)\,:\, 0\le t \le T\}$, and on the final time~$T$.
\end{theorem}

\begin{remark} Since the energy $H(x,v)=\tfrac12 |v|^2 + \phi(x)$ is conserved, it is bounded by $\tfrac12 \widetilde M^2 := \tfrac12 M^2 + \phi(x(0))$ and  we have $\tfrac12 |\dot x(t)|^2 \le  \widetilde M^2 - \phi(x(t))$.  As we assumed that the scalar potential $\phi$ is bounded from below, this gives an a priori bound on the velocity. Hence, the solution stays in a ball with centre $x(0)$ and radius depending only on $x(0)$ and $\dot x(0)$ in a fixed time interval $0\le t \le T$.
\end{remark}

\begin{remark} The differential equations for $z^0_0$ and $z^0_{\pm 1}$ are implicit, because the term $\ddot P(z^0)(z^0-c^0)$ contains $\ddot z^0_0$. By our choice of $c^0$, which ensures that $|z^0-c^0|$ is sufficiently small, the equation can be solved for $\ddot z^0_0$ to yield an explicit second-order differential equation. Similarly, the first-order differential equations for $z^0_{\pm1}$, which contain the time derivative in the term $\dot{P}_{\pm1}(z^0)(z^0-c^0)$, can be solved for $\dot z^0_{\pm 1}$ to yield explicit first-order differential equations. As was noted in \cite{hairer20lta}, the modulation functions $z^k$ are independent of the choice of~$c^0$.
\end{remark}

\begin{remark}
From the second equation of (c), it is straightforward to get (with $P_\perp=P_1+P_{-1}$ and $P_\parallel=P_0$)
\[
P_{\perp}\dot{z}^0=\frac{1}{|B|}P_{\parallel}\dot{z}^0\times \frac{\d b}{\d t}+\frac{1}{|B|^2}\left(E-\mu^0\,\nabla |B|\right)\times B+O(\eps^2),
\]
with $b=B/|B|, B, \nabla |B|, P_\parallel,P_\perp$ and $E$ evaluated at the guiding centre $z^0$,
which shows several slow drifts for the guiding centre motion usually derived by averaging techniques in the physical literature. 
\end{remark}

\begin{proof} It is sufficient to prove the theorem on time intervals of length $\eps$. At the end of an interval $[(n-1)\eps,n\eps]$, the construction of the modulated Fourier expansion is restarted from the exact solution values $x(n\eps),\dot x(n\eps)$, which in view of the uniqueness of the modulation functions up to $O(\eps^N)$ stated in (a) and the bound of the remainder term stated in (b) leads to jump discontinuities of size $O(\eps^N)$ in the modulation functions and the derivative of the phase function.

Statements (a) and (b) are given by Theorem 4.1 in~\cite{hairer20lta}. Here, we just give the proof of (c) and (d).

(c):
Inserting the modulated Fourier expansion into the differential equation~\eqref{ode} and comparing the coefficients of $\e^{\iu k\varphi(t)/\varepsilon}$ yields 
\[
\ddot{z}^k+2\iu k\frac{\dot{\varphi}}{\varepsilon}\dot{z}^k+\left(\iu k\frac{\ddot{\varphi}}{\varepsilon}-k^2\frac{\dot{\varphi}^2}{\varepsilon^2}\right)z^k=F^k,
\]
where the right-hand side $F^k$ is obtained from a Taylor expansion of $B$ and $E$ at $z^0$; see \cite{hairer20lta} for the general formula.
%\[
%F^k=\sum_{k_1+k_2=k}(\dot{z}^{k_1}+{\rm i}k_1\frac{\dot{\varphi}}{\varepsilon}z^{k_1})\times\sum_{m\geq0\atop s(\alpha)=k_2}\frac{1}{m!}B^{(m)}(z^0)\mathbf{z}^\alpha+\sum_{m\geq0\atop s(\alpha)=k}\frac{1}{m!}E^{(m)}(z^0)\mathbf{z}^\alpha.
%\]
For $k=0$, we obtain the motion of the guiding centre $z^0(t)$:
\begin{equation}\label{eq:z0}
\ddot{z}^0=\dot{z}^0\times B(z^0)+E(z^0)+\underbrace{2\,{\rm Re}\Bigl(\frac{{\rm i}\dot{\varphi}}{\varepsilon}\,z^{1}\times B'(z^{0})z^{-1}\Bigr)}_{=: I}%+2\,{\rm Re}\Bigl(\dot{z}^1\times B'(z^0)z^{-1}\Bigr)
+O(\eps).
\end{equation}
For $k=\pm 1$, we have
\begin{equation}\label{eq:z1}
\pm 2\iu\frac{\dot{\varphi}}{\varepsilon}\dot{z}^{\pm1}+\left(\pm\iu\frac{\ddot{\varphi}}{\eps}-\frac{\dot{\varphi}^2}{\eps^2}\right)z^{\pm1}=\left(\dot{z}^{\pm1}\pm \iu\frac{\dot{\varphi}}{\eps}z^{\pm1}\right)\times B(z^0)+\dot{z}^0\times B'(z^0)z^{\pm1}+O(\eps).
\end{equation}
%Comparing the leading terms on both hand sides of~\eqref{z1}, we have
%\[
%-\frac{\dot{\varphi}^2}{\eps^2}z^{\pm1}=\pm \iu\frac{\dot{\varphi}}{\eps^2}z^{\pm1}\times B(z^0)+O(1).
%\]

We first study the case $k=0$, i.e. \eqref{eq:z0}. Here we begin by giving an alternative expression for the term $I$, which is an $O(1)$ term. We show that
\begin{equation} \label{term-I}
I=-\mu^0\,\nabla|B|(z^0)+O(\eps).
\end{equation}
With the normalized eigenvectors $\nu_j$, we have $z_1^1=\zeta\nu_1$ and $z_{-1}^{-1}=\overline\zeta \nu_{-1}$ with $\nu_{-1}=\overline{\nu}_1$. We define the local orthonormal basis $e_1, e_2, e_3$ of $\mathbb{R}^3$ by the eigenvectors $\nu_j$ as $\nu_0= B/|B|=e_1$ and $\nu_{\pm 1}=\frac{1}{\sqrt 2}(e_2\pm \iu e_3)$. Using that $z^k_j=O(\eps^2)$ for $|k|=1$ and $j\ne k$ by part (a), the term $I$ can then be written as
\begin{align}\nonumber
I&=\frac{{\rm i}\dot{\varphi}}{\varepsilon}\,z^{1}_1\times  B'(z^{0})z^{-1}_{-1}-\frac{{\rm i}\dot{\varphi}}{\varepsilon}\,z^{-1}_{-1}\times  B'(z^{0})z^{1}_{1} + O(\eps)\\
\label{I-e}
&=|B(z_0)||z^1_1|^2 \left(e_2\times B'(z^0)e_3-e_3\times B'(z_0)e_2\right) + O(\eps).
\end{align}
Following equation (11) in~\cite{northrop63tam}, we find
\begin{align}\label{eq:northrop}
e_2\times B'(z^0)e_3-e_3\times B'(z^0)e_2=-\nabla|B|(z^0).
\end{align}
On the other hand, 
\[
x=z^0+O(\eps), \quad \dot{x}=\dot{z}^0+\iu\frac{\dot{\varphi}}{\eps}z^1_1\e^{\iu\varphi/\eps}-\iu\frac{\dot{\varphi}}{\eps}z^{-1}_{-1}\e^{-\iu\varphi/\eps}+O(\eps)
\]
and thus
\[
\dot{x}\times B(x)=\dot{z}^0\times B(z^0)-|B(z^0)|^2\left(z^1_1\e^{\iu\varphi/\eps}+z_{-1}^{-1}\e^{-\iu\varphi/\eps}\right)+O(\eps^0).
\]
From the orthogonality of $z_1^1$ and $z_{-1}^{-1}$ it follows that
\begin{align}\label{eq:moment}
\mu(x,\dot{x})=\frac{1}{2}\frac{|\dot{x}\times B(x)|^2}{|B(x)|^3}=|B(z^0)|\,{|z_1^1|^2}+O(\eps^2).
\end{align}
Inserting \eqref{eq:northrop} and \eqref{eq:moment} into \eqref{I-e} gives
\[
I=-\mu(x,\dot{x})\nabla|B|(z^0)+O(\eps).
\]
Using the adiabatic invariance \cite{hairer20lta,northrop63tam}
$
\: \mu(x(t),\dot x(t)) = \mu^0 + O(\eps^2),
$
we obtain \eqref{term-I}, and hence 
equation \eqref{eq:z0} can be equivalently written as
\begin{equation}\label{eq:z0new}
\ddot{z}^0=\dot{z}^0\times B(z^0)+E(z^0)-\mu^0\,\nabla|B|(z^0)+O(\eps).
\end{equation}
--- Multiplying~\eqref{eq:z0new} with $P_0(z^0)$ gives
\begin{align*}
P_0(z^{0})\ddot{z}^{0}=P_0(z^{0})\bigl(E(z^{0})-\mu^0\,\nabla |B|(z^{0})\bigr)+O(\eps).
\end{align*}
Using the product rule
$$\ddot z^0_0 = \frac{d^2}{dt^2} \bigl( P_0(z^0)(z^0-c^0) \bigr) =
P_0(z^{0})\ddot{z}^{0} + 2 \dot P_0(z^{0})\dot{z}^{0} + \ddot P_0(z^{0})({z}^{0}-c^0),
$$
this gives the first equation in (c).

\noindent
--- Multiplying~\eqref{eq:z0new} with $P_{\pm1}(z^0)$ gives 
\[
P_{\pm1}(z^0)\ddot{z}^0=\pm\iu\frac{\dot{\varphi}}{\eps}P_{\pm1}(z^0)\dot{z}^0+P_{\pm1}(z^0)(E(z^0)-\mu^0\,\nabla |B|(z^{0}))+O(\eps).
\]
Substituting $P_{\pm1}(z^0)\dot{z}^0=\dot{z}^0_{\pm1}-\dot{P}_{\pm1}(z^0)(z^0-c^0)$ yields
\begin{align*}
\dot{z}^0_{\pm1}-\dot{P}_{\pm1}(z^0)(z^0-c^0)&=\mp\iu\frac{\eps}{\dot{\varphi}}P_{\pm1}(z^0)\ddot{z}^0\pm\iu\frac{\eps}{\dot{\varphi}}P_{\pm1}(z^0)\left(E(z^0)-\mu^0\,\nabla |B|(z^{0})\right)+O(\eps^2)\\
&=\mp\iu\frac{\eps}{ \dot{\varphi}}\left(\ddot{z}^0_{\pm1}-\ddot{P}_{\pm1}(z^0)(z^0-c^0)-2\dot{P}_{\pm1}(z^0)\dot{z}^0\right)
\\
&\ \ \pm\iu\frac{\eps}{\dot{\varphi}}P_{\pm1}(z^0)\left(E(z^0)-\mu^0\,\nabla |B|(z^{0})\right)+O(\eps^2).
\end{align*}
Denoting $g_{\pm1}=\dot{z}^0_{\pm1}-\dot{P}_{\pm1}(z^0-c^0)$, we have $\dot{g}_{\pm1}=\ddot{z}^0_{\pm1}-\ddot{P}_{\pm1}(z^0)(z^0-c^0)-\dot{P}_{\pm1}(z^0)\dot{z}^0$. The above equation can be expressed as
\[
g_{\pm1}=\mp\iu\frac{\eps}{ \dot{\varphi}}\dot{g}_{\pm1}\pm\iu\frac{\eps}{\dot{\varphi}}\dot{P}_{\pm1}(z^0)\dot{z}^0\pm\iu\frac{\eps}{\dot{\varphi}}P_{\pm1}(z^0)\left(E(z^0)-\mu^0\,\nabla |B|(z^{0})\right)+O(\eps^2).
\]
By differentiation and substitution, the first term on the right-hand side can be absorbed into the $O(\eps^2)$ term, and so we get the second equation in (c).

Since the $\eps^{-2}$-terms cancel in~\eqref{eq:z1} after projection with $P_{\pm 1}(z^0)$, the $\eps^{-1}$-terms are dominant and we obtain the last equation in (c). 

(d): The initial values can be obtained by the same arguments as in the proof of Theorem 4.1 in~\cite{hairer20afb}.
\qed
\end{proof}

\subsection{Modulated Fourier expansion of the numerical solution}
The modulated Fourier expansion can be extended to the numerical solution of the modified Boris algorithm similarly to Theorem 4.2 in~\cite{hairer22lsi}. There are, however, additional terms and difficulties to be considered, since here we do not have a magnetic field in a near-constant direction as in~\cite{hairer22lsi}.
%We divide the interval $[0,T]$ into small intervals of length $O(h)$ and denote $t_m=mh$.
\begin{theorem}\label{thm:mfe-boris}
Let $x^n$ be the numerical solution obtained by applying the modified Boris algorithm to \eqref{ode}--\eqref{xv-bounds} with a stepsize $h$ satisfying
\begin{equation}\label{hh-eps}
c_*\eps \le h^2 \le C_*\eps
\end{equation}
for some positive constants $c_*$ and $C_*$.
We assume that the component orthogonal to $B(x^0)$ 
of the starting velocity, $v^0_\perp=P_\perp(x^0) v^0= v^0 - P_0(x^0)v^0$, is chosen to be small:
\begin{equation}\label{v-init-boris}
|v^0_\perp| \le c_1 \eps .
\end{equation}
We further make the nondegeneracy assumption \eqref{ass-nondeg}.
%\begin{align}\label{ass-nondeg}
%&\text{For $(x,v)$ along the numerical trajectory, the linear maps} \nonumber
%\\
%& \text{$L_{x,v}:P_\perp(x)\mathbb{R}^3 \to P_\perp(x)\mathbb{R}^3, \quad z \mapsto z + \tfrac14 h^2\, P_\perp(x)\bigl(v \times B'(x)z\bigr)$}
%\\
%&\text{are invertible with an inverse bounded by a constant.}  \nonumber
%\end{align}
%We further assume that the numerical solution $x^n$ stays in a compact set $K$ for $0\leq nh\leq T_{h}$ with $T_h=O(h)$. 
For an arbitrary truncation index $N\geq 2$, we then have a decomposition
\begin{equation} \label{mfe-boris}
x^n=y(t_n) + (-1)^n z(t_n) +R_N(t_n), \qquad t_n=nh \le T,
\end{equation}
%with {\blue $y(t)=\sum_m y^{[m]}(t)\mathbbm 1_m(t)$ and $z(t)=\sum_m z^{[m]}(t)\mathbbm 1_m(t)$.}
with the following properties:
\begin{itemize}
\item [(a)] 
The  functions $y(t)$ and $z(t)$, $0\le t \le T$, are piecewise continuous with jumps of size $O(h^N)$ at integral multiples of $h$ and are smooth elsewhere. Together with their derivatives (up to order $N$) they are bounded as $y=O(1)$, $z=O(h^2)$. They are unique up to $O(h^{N})$. 
Moreover, $P_\perp(y)\dot y = O(\eps)$ and $P_0(y)z=O(h^4)$.
%$\dot{y}\times B(y)=O(\varepsilon^0)$ and {\red $P_0(y)(\dot{y}\times B'(y)z) =O(h^2)$}.
\item [(b)] 
The remainder term is
%and its derivative are 
bounded by
\[
R_N(t)=O(h^N) %,\quad \dot{R}_N(t)=O(t h^N) 
\quad \text{for} \quad 0\leq t\leq T.
\]
\item [(c)]  We let $c^0(t)$ be a piecewise constant function that is sufficiently close to $y(t)$.
The functions $y_j=P_j(y) (y- c^0)$ $(j=0,\pm1)$ and $z_{\pm 1}=P_{\pm1}z$ satisfy the following differential equations for $0\le t \le T$ except at the jumps. Here, all functions $B$, $E$, $P_j$ are evaluated at the numerical guiding centre $y(t)$, and we write $\dot P_j = (d/dt) P_j(y(t))=P_j'(y(t)) \dot y(t)$ and analogously $\ddot P_j$. Moreover, $\mu^0=\mu(x(0),\dot x(0)))$ is the magnetic moment. Omitting the ubiquitous argument $t$, we have
%\[
%\begin{aligned}
%P_0(y^{[m]})\ddot{y}^{[m]}=&P_0(y^{[m]})\red{\overbrace{\bigl(-\dot{z}^{[m]}\times \frac{1}{\eps}B'(y^{[m]})z^{[m]}\bigr)}^{=\mu\nabla|B(y^{[m]})|=O(\eps^2)}}\\
%&+P_0(y^{[m]})E(y^{[m]})-P_0(y^{[m]})\bigl(\bar{\mu}\nabla|B(y^{[m]})|\bigr)+O(h^2),\\
%P_{\pm 1}(y^{[m]})\dot{y}^{[m]}=&\pm{\rm i}\frac{\eps}{|B(y^{[m]})|} P_{\pm 1}(y^{[m]})\red{\bigl(-\dot{z}^{[m]}\times \frac{1}{\eps}B'(y^{[m]})z^{[m]}\bigr)}\\
%&\pm{\rm i}\frac{\eps}{|B(y^{[m]})|} P_{\pm 1}(y^{[m]})\bigl(E(y^{[m]})-\bar{\mu}\nabla|B(y^{[m]})|\bigr)+O(\eps h^2),\\[1mm]
%P_{\pm 1}(y^{[m]})\dot{z}^{[m]}=&\mp 4\iu \frac\eps{h^2|B(y^{[m]})|} {z}^{[m]}_{\pm 1} \mp \iu \frac{1}{|B(y^{[m]}|}P_{\pm 1}(y^{[m]})\bigl(\dot{y}^{[m]}\times B'(y^{[m]})\bigr)+ O(\eps h^2).
%\end{aligned}
%\]
%\begin{align*}
%\ddot{y}_0&=P_0(y)\left(E(y)-\mu^0\,\nabla|B|(y)\right)+2\dot{P}_0(y)\dot{y}+\ddot{P}_0(y)(y-c^0)+O(h^2)\\
%\dot{y}_{\pm1}&=\dot{P}_{\pm1}(y)(y-c^0) \pm\frac{\iu}{|B(y)|}\dot{P}_{\pm1}(y)\dot{y}
%\\
%&\ \ 
%\pm \frac{\iu}{|B(y)|}P_{\pm1}(y)\left(E(y)-\mu^0\,\nabla|B|(y)\right)+O(\eps h^2)\\
%\dot{z}_{\pm1}&=\dot{P}_{\pm1}(y)z\mp\frac{4\iu}{h^2|B(y)|}z_{\pm1}{\red \mp\,\frac{\iu}{|B(y)|}P_{\pm1}(y)\left(\dot{y}\times B'(y)z\right)}+O(\eps h^2).
%\end{align*}
\begin{align*}
\ddot{y}_0&=2\dot{P}_0\dot{y}+\ddot{P}_0(y-c^0)+ P_0\left(E-\mu^0\,\nabla|B|\right)+O(h^2)\\
\dot{y}_{\pm1}&=\dot{P}_{\pm1}(y-c^0) \pm\frac{\iu}{|B|}\dot{P}_{\pm1}\dot{y}
\pm \frac{\iu}{|B|}P_{\pm1}\left(E-\mu^0\,\nabla|B|\right)+O(h^2)\\
\dot{z}_{\pm1}&=\dot{P}_{\pm1}z\mp\frac{4\iu}{h^2|B|}z_{\pm1}{\mp\,\frac{\iu}{|B|}P_{\pm1}\left(\dot{y}\times B'(y)z\right)}+O(\eps h^2).
\end{align*}

The function $z_0=P_0(y)(z-c_0)$ is given by an algebraic expression depending on $y$, $\dot{y}_0$ and $z_{\pm 1}$.
\item [(d)]  %With the choice $c^0(0)=x^0$,
Initial values for the differential equations of item (c) are given by
\[
\begin{aligned}
y(0)&= x^0+O(h^2),\\
\dot{y}_0(0)&=P_0(x^0) v^0 +O(h^2),
\\
z_{\pm 1}(0)&= O(h^2).
%\mp\frac{\iu}{4} h^2|B(x^0)| \,P_{\pm 1}(x^0)\Bigl( v^0 \mp \frac{\iu}{|B(x^0)|} \bigl(E(x^0)-{\mu^0}\nabla|B|(x^0)\bigr)\Bigr)+O(h^2).
\end{aligned}
\]
\end{itemize}
The constants symbolized by the $O$-notation are independent of $\varepsilon$, $h$ and $n$ with $0\leq nh \leq T$, but depend on the velocity bound, on bounds of derivatives of $B$ and $E$ in a neighbourhood of the numerical trajectory, and on the final time~$T$.
\end{theorem}

\begin{remark} The essential observation is that for the modified Boris method, the differential equations for  the numerical guiding centre $y(t)$ are the same, up to a defect of size $O(h^2)$, as the differential equations for  the  guiding centre $z^0(t)$ of the exact solution, and also the initial values agree up to $O(h^2)$. In contrast, for the standard Boris method with parallel-projected initial velocity, the terms $\mu^0\,\nabla|B|$ are missing. This is the reason for the failure of the standard Boris method with modified starting values for large step sizes $h^2\ge \eps$ in the situation of strongly non-uniform strong magnetic fields.
\end{remark}

\begin{proof}
This theorem is proved similarly to Theorem 4.2 in~\cite{hairer22lsi} (which gives an analogous decomposition for the standard Boris method in the case of a near-constant strong magnetic field)
combined with the treatment of the strongly nonuniform magnetic field in Theorem 4.1 in \cite{hairer20lta}. Here, we do not repeat the arguments in the proofs of those papers for (a) and (b) (such as the recursive elimination of higher time derivatives, an idea going back in time as far as the Euler--Maclaurin summation formula \cite{hairer97abi}) but concentrate on the parts (c) and (d) that are specific for the present situation.

Since a general strong magnetic field is considered, the time interval of validity of the modulated Fourier expansion is here $O(h)$ instead of $O(1)$, and so we need to patch together many such short-time expansions, starting anew from each $x^n$, in the same way we did in Theorem~\ref{thm:mfe} over intervals of length proportional to $\eps$.
%Over a time interval of length $O(h)$, an application of discrete Gronwall inequality allows us to conclude from a defect of size $O(h^N)$ to an error of size $O(t^2h^N)$.

Inserting the decomposition~\eqref{mfe-boris} into the numerical method~\eqref{mboris} and separating the terms without and with the factor $(-1)^n$ gives
\begin{align}
\label{eq:y}
&\ddot{y}+O(h^2)=\bigl(\dot{y}+ O(h^2)\bigr)\times B(y)-\dot{z}\times B'(y)z+E(y)-\mu^0\,\nabla|B|(y)+O(h^2)\\
\label{eq:z}
&-\frac{4}{h^2}z-\ddot{z}+O(h^2)=-\dot{z}\times B(y)+\dot{y}\times B'(y)z+E'(y)z+O(h^2).
\end{align}
%As in the proof of Theorem~\ref{thm:mfe}, the 
Since $z=O(h^2)$ and $\dot z=O(h^2)$, the second term on the right hand side of the first equation is
\[
\dot{z}\times B'(y)z=  O(h^4/\eps) =   O(h^2)   
\]
in our stepsize regime $h^2\sim\eps$.

Taking the projection $P_0=P_0(y)$ on both sides of~\eqref{eq:y} yields the first equation in (c). Taking the projection $P_{\pm1}$ on both sides gives
\[
P_{\pm1}\ddot{y}+O(h^2)=\pm\iu |B(y)| P_{\pm1}\dot{y}+P_{\pm1}\left(E(y)-\mu^0\,\nabla|B|(y)\right)+O(h^2 |B(y)|).
\]
As in Theorem~\ref{thm:mfe} we thus have, with $B=B(y)$,
\[
P_{\pm1}\dot{y}=\pm\frac{\iu}{|B|}\dot{P}_{\pm1}\dot{y}\pm\frac{\iu}{|B|}P_{\pm1}\left(E(y)-\mu^0\,\nabla|B(y)|\right)+O( h^2).
\]

Taking the projection $P_{\pm1}=P_{\pm 1}(y)$ on both sides of~\eqref{eq:z} yields
\[
-\frac{4}{h^2}z_{\pm1}-P_{\pm1}\ddot{z}+O(h^2)=\mp\iu |B|P_{\pm1}\dot{z}+P_{\pm1}\left(\dot{y}\times B'(y)z\right)+O(h^2),
\] 
and so we find
\[
P_{\pm1}\dot{z}=\mp\frac{4\iu}{h^2|B|}z_{\pm1}\mp\frac{\iu}{|B|}P_{\pm1}\left(\dot{y}\times B'(y)z\right)+O(\eps h^2).
\]
We thus have the differential equations of part (c). Taking $P_0$ on both sides of~\eqref{eq:z} and multiplying with $-h^2/4$ yields
$$
z_0=-\tfrac14 h^2\,P_0(\dot{y}\times B'(y)z)+O(h^4).
$$
Since $P_\perp \dot y=O(\eps)$, we have $P_0(\dot{y}\times B'(y)z)= P_0(P_\perp\dot{y}\times B'(y)z)=O(z)$. This 
gives us $z_0=O(h^4)$ provided that $z_{\pm1}=O(h^2)$.

(d) The numerical approximation to the velocity is given by
\[
v^n=\frac{x^{n+1}-x^{n-1}}{2h}=\dot{y}(t_n)+\dddot{y}(t_n)h^2+\cdots-(-1)^n(\dot{z}(t_n)+\dddot{z}(t_n)h^2+\cdots),
\]
and so we have
\[
v^n_{\perp}=P_{\perp}\dot{y}(t_n)-(-1)^n P_{\perp}\dot{z}(t_n)+ O(h^2),
\]
which under the bounds of (a) yields $v^n_{\perp}=O(h^2)$. We now consider this equation for $n=0$.
Since the above equation for $P_{\pm1}\dot{z}$ and the bound for $z_0$ yield
$$
P_\perp \dot z(0) = \frac{4}{h^2|B^0|}\, L_{x^0,v^0}(z_\perp(0))  \times \frac{B^0}{|B^0|} + O(h^4),
$$
the above equation for $v^0_\perp$ yields
$$
\frac4{h^2|B^0|}\, L_{x^0,v^0}(z_\perp(0)) \times \frac{B^0}{|B^0|}= P_{\perp}\dot{y}(0) - v^0_\perp + h^2 P_{\perp}\dddot{y}(0) + O(h^2 z) +O(h^4),
$$
and with the nondegeneracy condition \eqref{ass-nondeg} we are now able to construct $z_\perp(0)$ and hence $z_{\pm 1}(0)$, which thanks to $h^2\sim\eps$ and $v^0_\perp=O(\eps)$ are indeed of size $O(h^2)$.
% [We remark that against expectations, the most problematic term for an extension of the result to step sizes $h^2 \gg\eps$ is the term $h^2 P_{\perp}\dddot{y}(0)$, which is only $O(h^2)$.]
\qed
\end{proof}

\subsection{Proof of Theorem 2.1}

Theorem~\ref{thm:mfe} represents the exact solution as
$$
x(t)= z^0(t) + O(\eps),
$$
and Theorem~\ref{thm:mfe-boris} represents the numerical solution of the modified Boris method with $h^2\sim\eps$ as
$$
x^n = y(t_n) + O(h^2),
$$
where the guiding centre $z^0(t)$ and the numerical guiding centre $y(t_n)$ satisfy the same differential equations up to $O(h^2)$ with the same initial values up to $O(h^2)$, and the jumps of size $O(\eps^N)$ or $O(h^N)$ for arbitrary $N$ contribute less than $O(h^2)$ to the difference. (The piecewise constant function $c^0(t)$ can be chosen the same in both cases.) Therefore, $z^0(t)$ and $y^0(t)$ differ by $O(h^2)$ on a fixed time interval $0\le t \le T$. This proves the $O(h^2)$ error bound for the positions in Theorem~\ref{thm:main}.

We now turn to the error bound for the velocity. We compare the velocity of the exact solution 
\[
v(t)=\dot{x}(t)=\dot{z}^0(t)+\frac{\iu \dot{\varphi}(t)}{\eps}z^1_1(t)\e^{\iu\varphi(t)/\eps}-\frac{\iu \dot{\varphi}(t)}{\eps}z^{-1}_{-1}(t)\e^{-\iu\varphi(t)/\eps}+O(\eps)
\]
and the numerical velocity 
\[
v^n=\frac{x^{n+1}-x^{n-1}}{2h}=\dot{y}(t_n)-(-1)^n\dot{z}(t_n)+O(h^N).
\]
Since $P_\parallel(z^0)z^{\pm1}_{pm1}=0$ and $P_\parallel(y)z=z_0=O(h^4)$, and since we already know that $z^0(t)-y(t)=O(h^2)$ and
$z^0(t)-x(t)=O(\eps)$ and $y(t_n)-x^n=O(h^2)$,
it follows that
\[
v_{\parallel}^n-v_{\parallel}(t^n)= P_\parallel(x^n)v^n - P_\parallel(x(t_n))v(t_n) = O(h^2).
\]
Finally, the bound $v^n_\perp=O(h^2)$ was already shown in part (d) of the proof of Theorem~\ref{mfe-boris}.
This completes the proof of Theorem~\ref{thm:main}.

\section*{Acknowledgement} 
We thank Michael Kraus for making us aware of Xiao and Qin's modified Boris integrator in \cite{xiao21smc} and the slow manifold approach in \cite{burby21nso}.
This work was supported by the Deutsche Forschungsgemeinschaft (DFG, German Research Foundation) -- Project-ID258734477 -- SFB 1173 and a joint DAAD--CSC postdoctoral scholarship.

\bibliographystyle{abbrv}
\bibliography{HLW}

\end{document}